\documentclass[9pt,journal]{IEEEtran}
\usepackage{amssymb,latexsym,amsmath,cite,array,rsfso,graphicx}
\usepackage[ampersand]{easylist}
\newtheorem{prop}{Proposition}[section]
\newtheorem{dfn}{Definition}[section]
\newtheorem{exm}{Example}[section]
\newtheorem{rem}{Remark}[section]
\numberwithin{equation}{section}

\graphicspath{{Observability and State Estimation for a Class of Nonlinear Systems - FINAL arXiv}}

\begin{document}
\title{Observability and State Estimation for a Class of Nonlinear Systems}
\author{\IEEEauthorblockN{John Tsinias and
Constantinos Kitsos}

\thanks{J. Tsinias is with the Department of Mathematics, National Technical
University of Athens, Zografou Campus 15780, Athens, Greece, email:{\tt\small jtsin@central.ntua.gr} (corresponding author). C. Kitsos is with the Department of Mathematics, National Technical
University of Athens, Zografou Campus 15780, Athens, Greece and Univ. Grenoble Alpes, CNRS, Grenoble INP$^\circ$, GIPSA-lab, 38000 Grenoble, France, email: {\tt\small konstantinos.kitsos@gipsa-lab.grenoble-inp.fr}

{\tiny $\circ$ Institute of Engineering Univ. Grenoble Alpes}}}

\markboth{}%
{Shell \MakeLowercase{\textit{et al.}}: Bare Demo of IEEEtran.cls for IEEE Journals}

\maketitle

\begin{abstract}
We derive sufficient conditions for the solvability of the state estimation problem for a class of nonlinear control time-varying systems which includes those, whose dynamics have triangular structure. The state estimation is exhibited by means of a sequence of functionals approximating the unknown state of the system on a given bounded time interval. More assumptions guarantee solvability of the state estimation problem by means of a hybrid observer.  
\end{abstract}


\IEEEpeerreviewmaketitle

\section{Introduction}
Many important approaches have been presented in the literature concerning the state estimation problem for a given nonlinear control system (see for instance [1] - [19], [21] - [25], [27], [28] and relative references therein). Most of them are based on the existence of an observer system exhibiting state estimation. The corresponding hypotheses include observability assumptions and persistence of excitation. In [8], [9] and [10] Luenberger type observers and switching estimators are  proposed for a general class of triangular systems under weaker assumptions than those adopted in the existing literature. In [16] and [24] the state estimation is exhibited for a class of systems by means of a hybrid observer.

The present note is inspired by the approach adopted in [16], where a "hybrid dead-beat observer" is used, as well as by the methodologies applied in [20], [26] and [29], where, fixed point theorems are used for the establishment of sufficient conditions for observability and asymptotic controllability. Our main purpose is to establish that, under certain hypotheses, including persistence of excitation, the State Estimation Design Problem (SEDP) around a given fixed value of initial time is solvable for a class of nonlinear systems by means of a \textit{sequence of mappings} $X_{\nu}$,  exclusively dependent on the dynamics of the original system, the  input $u$ and the corresponding output $y$, and further each $X_{\nu}$ is independent of the time-derivatives of $u$ and $y$. An algorithm for explicit construction of these mappings is provided.

We consider time-varying finite dimensional nonlinear control systems of the form:
\begin{subequations}
\begin{IEEEeqnarray}{rCl}       
\begin{aligned}       
       \dot{x}=A(t,y,u)x+f(t,y,x,u)\IEEEeqnarraynumspace\\ 
        (t,y,x,u)\in {{\mathbb{R}}_{\ge 0}}\times {{\mathbb{R}}^{k}}\times {{\mathbb{R}}^{n}}\times {{\mathbb{R}}^{m}}
\end{aligned}\IEEEeqnarraynumspace\IEEEeqnarraynumspace\IEEEeqnarraynumspace\\        
       \text{with output}\IEEEeqnarraynumspace\IEEEeqnarraynumspace\IEEEeqnarraynumspace\IEEEeqnarraynumspace\IEEEeqnarraynumspace\IEEEeqnarraynumspace\IEEEeqnarraynumspace\IEEEeqnarraynumspace\IEEEeqnarraynumspace\nonumber\\
       y=C\left( t,u\right) x \IEEEeqnarraynumspace\IEEEeqnarraynumspace \IEEEeqnarraynumspace\IEEEeqnarraynumspace \IEEEeqnarraynumspace    
\end{IEEEeqnarray}
\end{subequations}
where $u$ is the input of (1.1). Our main results  establish  sufficient conditions for the approximate solvability of the SEDP for (1.1). The paper is organized as follows. Section II contains definitions, assumptions, as well as statement and proof of our main result (Proposition 2.1) concerning the state estimation problem for the general case (1.1). We apply in Section III the main result of Section II for the derivation of sufficient conditions for the solvability of the same problem for certain  subclasses of systems (1.1), whose dynamics have triangular structure (Proposition 3.1). According to our knowledge, the sufficient conditions proposed in Sections II, III are weaker than those imposed in the existing literature for the solvability of the observer design problem for the same class of systems. More extensions are discussed in Section IV of present work, concerning the solvability of the SEDP by means of a hybrid observer under certain additional assumptions (Proposition 4.1).

\textit{Notations:} For given $x\in {{\mathbb{R}}^{n}}$, $\left| x \right|\text{ }$denotes its usual Euclidean norm. For a given constant matrix $A\in {{\mathbb{R}}^{m\times n}}$, ${A}'$ denotes its transpose and $\left| A \right|:=\sup \left\{ \left| Ax \right|,\left| x \right|=1 \right\}$ is its induced norm. For any nonempty set $I\subset \mathbb{R}$ and map $I\ni  t\to A(t)\in {{\mathbb{R}}^{m\times n}}$ we adopt the notation $\|A(\cdot ){{\|}_{I}}\text{:=}\operatorname{ess}.\sup\{\left| A(t) \right|,t\in I\}$.

\section{Hypotheses and Main Result}
 In this section we provide sufficient conditions for observability of (1.1) and solvability of the SEDP. We assume that for each $t\ge 0$ the mappings $A\left( t,\cdot ,\cdot  \right):$ ${{\mathbb{R}}^{k}}\times {{\mathbb{R}}^{m}}\to {{\mathbb{R}}^{n\times n}}$ , $C\left( t,\cdot  \right):{{\mathbb{R}}^{m}}\to {{\mathbb{R}}^{k\times n}}$ and $f\left( t,\cdot ,\cdot ,\cdot  \right):$ ${{\mathbb{R}}^{k}}\times {{\mathbb{R}}^{n}}\times {{\mathbb{R}}^{m}}\to {{\mathbb{R}}^{n}}$ are continuous and further $f$ is (locally) Lipschitz continuous with respect to $x$, i.e., for every bounded $I\subset {{\mathbb{R}}_{\ge 0}},X\subset {{\mathbb{R}}^{n}},U\subset {{\mathbb{R}}^{m}}$ and $Y\subset {{\mathbb{R}}^{k}}$ there exists a constant $C>0$ such that 
\begin{IEEEeqnarray}{rCl}
\left| f\left( t,y,{{z}_{1}},u \right)-f\left( t,y,{{z}_{2}},u \right) \right|\le C\left| {{z}_{1}}-{{z}_{2}} \right|,\nonumber \\
\forall \left( t,y,u \right)\in I\times Y\times U,\,{{z}_{1}},{{z}_{2}}\in X
\end{IEEEeqnarray}
Also, assume that for any $\left( x,y,u \right)\in {{\mathbb{R}}^{n}}\times {{\mathbb{R}}^{k}}\times {{\mathbb{R}}^{m}}$ the mappings $A\left( \cdot,y,u \right)$, $f\left( \cdot ,y,x,u \right)$ and $C(\cdot ,u)$ are measurable and locally essentially bounded in $\mathbb{R}_{\ge 0}$. Let ${{t}_{0}}\ge 0,\tau >{{t}_{0}}$ and let ${{U}_{[{{t}_{0}},\tau ]}}$ be a nonempty set of inputs $u\in {{L}_{\infty }}\left( \left[ {{t}_{0}},\tau  \right];{{\mathbb{R}}^{m}} \right)$ of (1.1) (without any loss of generality it is assumed that ${{U}_{[{{t}_{0}},\tau ]}}$ is independent of the initial state). Define by ${{Y}_{[{{t}_{0}},\tau ],u}}$ the set of outputs $y$ of (1.1) defined on the interval $[{{t}_{0}},\tau ]$ corresponding to some input $u\in {{U}_{[{{t}_{0}},\tau ]}}$:
\begin{IEEEeqnarray}{rCl}
Y_{[t_{0},\tau],u} :=\lbrace y \in L_{\infty}\left([t_{0},\tau];\mathbb{R}^{k} \right): y(t)=C\left(t,u(t)\right) \nonumber \\\times x\left(t;t_{0},x_{0},u\right), \text{a.e. } t \in [t_0,\tau], \text{ for certain } x_{0}\in\mathbb{R}^{n}\rbrace \IEEEeqnarraynumspace
\end{IEEEeqnarray}
provided that ${{t}_{\max }}\ge \tau $ where $t_{\max}=t_{\max}(t_{0},x_{0},u)\le + \infty$ is the maximum time of existence of the solution $x(\cdot ;{{t}_{0}},{{x}_{0}},u)$ of (1.1) with initial $x({{t}_{0}};{{t}_{0}},{{x}_{0}},u)={{x}_{0}}$.\\
\begin{dfn}
Let $I$ be a nonempty subset of ${{\mathbb{R}}_{\ge 0}}$. We say that  (1.1) is  \textit{observable over}  $I$, if for all ${{t}_{0}}\in I$, almost all $\tau >{{t}_{0}}$ near ${{t}_{0}}$, input  $u\in {{U}_{[{{t}_{0}},\tau ]}}$ and  output $y\in {{Y}_{[{{t}_{0}},\tau ],u}}$, there exists a unique ${{x}_{0}}\in {{\mathbb{R}}^{n}}$  such that 
\begin{IEEEeqnarray}{rCl}
y(t)=C(t,u(t))x(t;{{t}_{0}},{{x}_{0}},u),\text{a.e. } t\in \left[ {{t}_{0}},\tau  \right]\IEEEeqnarraynumspace
\end{IEEEeqnarray}
\end{dfn}
\par According to Definition 2.1, observability is equivalent to the existence of a (probably noncausal) functional $X(\cdot ,\cdot ,\cdot ):\left\{ {{t}_{0}} \right\}\times {{Y}_{[{{t}_{0}},\tau ],u}}\times {{U}_{[{{t}_{0}},\tau ]}}\to {{\mathbb{R}}^{n}}$, such that, for every ${{x}_{0}}\in {{\mathbb{R}}^{n}}$ for which (2.3) holds for certain $u\in {{U}_{[{{t}_{0}},\tau ]}}$ and $y\in {{Y}_{[{{t}_{0}},\tau ],u}}$, we have: 
\begin{IEEEeqnarray}{rCl}
X({{t}_{0}},y,u)={{x}_{0}} \IEEEeqnarraynumspace
\end{IEEEeqnarray}
and $X$ is exclusively dependent on the input $u$ and the output $y$ of (1.1) and is in general \textit{noncausal}. Knowledge of $X$ satisfying the previous properties guarantees knowledge of the initial state value, thus knowledge of the future values of the solution of (1.1), provided that the system is complete. 
\begin{dfn}
 Let $I$ be a nonempty subset of ${{\mathbb{R}}_{\ge 0}}$. We say that  the SEDP is  solvable for (1.1) over  $I$, if there exists a functional $X(\cdot,\cdot,\cdot ):\left\{ {{t}_{0}} \right\}\times {{Y}_{[{{t}_{0}},\tau ],u}}\times {{U}_{[{{t}_{0}},\tau ]}}\to {{\mathbb{R}}^{n}}$, ${{t}_{0}}\in I$, $\tau >{{t}_{0}}$ near ${{t}_{0}}$, being in general noncausal, such that (2.4) is fulfilled for every ${{x}_{0}}\in {{\mathbb{R}}^{n}}$ for which (2.3) holds for certain $u\in {{U}_{[{{t}_{0}},\tau ]}}$ and further $X$ \textit{is exclusively depended on $u$ and $y$ and the dynamics of (1.1) and it does not include any differentiation of their arguments}. It turns out that $X$ is independent of  the time-derivatives of $u$ and $y$, whenever they exist. 
\end{dfn}

It is worthwhile to remark here that the approach proposed in [16] for the construction of a hybrid dead-beat observer for a subclass of systems (1.1) is based on an explicit construction of a (noncausal) map $X$ satisfying (2.4). However, for general nonlinear systems, the precise and direct determination of the functional $X$ is a difficult task. The difficulty comes from our requirements for the candidate $X$ to be exclusively  dependent on $u$ and $y$ and the dynamics of system and, for practical reasons, it should be independent of the time-derivatives of $u$ and $y$. We next provide a weaker  sequential type of definition of the solvability of SEDP, which is adopted in the present work, in order to achieve the state determination for general case (1.1) by employing an explicit approximate strategy. 
\begin{dfn}
We say that  the  \textit{approximate SEDP}  is solvable for system (1.1) over $I$, if there exist functionals $X_\nu (\cdot ,\cdot ,\cdot ), \nu=1,2,\ldots:\left\{ {{t}_{0}} \right\}\times {{Y}_{[{{t}_{0}},\tau ],u}}\times {{U}_{[{{t}_{0}},\tau ]}}\to {{\mathbb{R}}^{n}}$(being in general noncausal), such that, if we denote:
\begin{IEEEeqnarray}{rCl}
{{\xi }_{\nu}}:={{X}_{\nu}}({{t}_{0}},y,u),\ {{t}_{0}}\in I, y\in {{Y}_{[{{t}_{0}},\tau ],u}}, u\in {{U}_{[{{t}_{0}},\tau ]}} 
\end{IEEEeqnarray}
then\\
(I) the mappings $X_{\nu}$ are exclusively dependent on the input $u$ and the corresponding output $y$, the dynamics of system (1.1) and further their domains do not include any differentiation of their arguments. It turns out that each $X_{\nu}$ should be  independent of  the time-derivatives of $u$ and $y$ (whenever they exist);\\
(II) the following hold: 
\begin{IEEEeqnarray}{CC}
\IEEEyesnumber \IEEEyessubnumber*
       \lim_{\nu \to \infty} {{\xi }_{\nu}}={{x}_{0}};\IEEEeqnarraynumspace\\
\begin{aligned}       
       {{x}_{0}} \text{ is the (unique) vector for}\IEEEeqnarraynumspace\IEEEeqnarraynumspace\\       \text{which both (2.3) and (2.6a) hold} \IEEEeqnarraynumspace
       \end{aligned}
       \end{IEEEeqnarray}
\end{dfn}
\par It should be emphasised that uniqueness requirement in (2.6b) is not essential. We may replace  (2.6b) by the assumption that there exists $x_0$ satisfying both (2.3) and (2.6a). Then uniqueness of such a vector $x_0$ is a consequence  of (2.6a), definition (2.5) and the fact that each functional $X_{\nu}$ exclusively depends on $u$ and $y$.
\par Obviously, according to the definitions above, the following implications hold:
\par Solvability of SEDP $\Rightarrow$ Solvability of approximate SEDP $\Rightarrow$ Observability (over $I$).\\
For completeness, we note that the first implication follows by setting $X_{\nu}:=X,\nu=1,2,\ldots$ in (2.4). The second implication is a direct consequence of both assumptions (2.6a,b), definition (2.5) and the exclusive dependence of each $X_{\nu}$ from $u$ and $y$. The converse claims are not in general valid; particularly, observability does not in general imply solvability of the (approximate) SEDP, due to the additional requirements of Definitions 2.2 and 2.3 concerning the independence $X$, $X_{\nu}$, respectively, from the time-derivatives of $u$ and $y$.
\par From (2.6a) we deduce that, if the approximate SEDP is solvable for (1.1) over $I$, then for any $T>t_0$ for which $T\le {{t}_{\max }}({{t}_{0}},{{x}_{0}},u)$ it holds: 
\begin{IEEEeqnarray}{rCl}
\lim_{\nu\to\infty} \Vert x(\cdot;t_0,\xi_{\nu},u)-x(\cdot;t_0,x(t_0),u)\Vert_{[t_0,T]}=0\IEEEeqnarraynumspace
\end{IEEEeqnarray}  
\begin{rem} (i) Condition (2.7) guarantees that for any given interval $[{{t}_{0}},T]$ with $T\le {{t}_{\max }}({{t}_{0}},{{x}_{0}},u)$, the unknown solution $x(s;{{t}_{0}},x({{t}_{0}}),u),s\in [{{t}_{0}},T]$ of (1.1) is uniformly approximated by a sequence of  trajectories $\hat{x}$ of the system 
\begin{IEEEeqnarray}{rCl}\dot{\hat{x}}\left( t \right)=A(t,y,u)\hat{x}+f(t,y,\hat{x},u), \hat{x}({{t}_{o}})={{\xi }_{\nu}}, \nu=1,2,\ldots  \nonumber
\end{IEEEeqnarray}
with ${{\xi }_{\nu}},\nu=1,2,\ldots$ as given in Definition 2.2.\\ 
(ii) If  the  system (1.1) is \textit{complete}, then (2.7)  implies solvability of the approximate SEDP, thus, observability for  (1.1) over ${{\mathbb{R}}_{\ge 0}}$. Indeed, let ${{t}_{0}}\in I$ and without loss of generality consider arbitrary $s>{{t}_{0}}$. It follows by invoking the forward completeness assumption and (2.7) that $\lim_{\nu \to \infty} {{\hat{\xi }}_{\nu}}=x(s)$, where  ${{\hat{\xi }}_{\nu}}:=x(s;{{t}_{0}},{{\xi }_{\nu}},u),\nu=1,2,\ldots$ and simultaneously (2.3) and (2.6a) hold with $s$ and $x(s)$, instead of ${{t}_{0}}$ and ${{x}_{0}}=x({{t}_{0}})$, respectively. Moreover, due to the backward completeness, $x(s)$ is the unique vector for which $\lim_{\nu \to \infty} {{\hat{\xi }}_{\nu}}=x(s)$.  
\end{rem}

In order to state and establish our main result, we first require the following notations and additional assumptions for the dynamics of (1.1). Consider ${{t}_{0}}\in I$, $\tau >{{t}_{0}}$, $u\in {{U}_{[{{t}_{0}},\tau ]}}$, $y\in {{Y}_{[{{t}_{0}},\tau ],u}}$ and $d\in {{C}^{0}}\left( \left[ {{t}_{0}},\tau  \right];{{\mathbb{R}}^{m}} \right)$. We denote by $\Phi (t,{{t}_{0}})$ the fundamental matrix solution of
\begin{IEEEeqnarray}{CC}
\IEEEyesnumber \IEEEyessubnumber*
\frac{\partial }{\partial t}\Phi (t,{{t}_{0}})=A(t,y(t),u(t))\Phi (t,{{t}_{0}})\IEEEeqnarraynumspace\\
\Phi \left( {{t}_{0}},{{t}_{0}} \right)={{I}_{n\times n}} \IEEEeqnarraynumspace
\end{IEEEeqnarray}
and define the mappings:  
\begin{multline}
\Psi \left( t;{{t}_{0}},y,u \right):=\\\int\limits_{{{t}_{0}}}^{t}{\Phi '\left( s,{{t}_{0}} \right)C'\left( s,u\left( s \right) \right)C\left( s,u\left( s \right) \right)\Phi \left( s,{{t}_{0}} \right)ds}, t\in \left[ {{t}_{0}},\tau  \right] 
\end{multline} 
\begin{multline}
{{\Xi }}(t;{{t}_{0}},y,d,u):=\int_{{{t}_{0}}}^{t}{{{\Phi }'}}(\rho ,{{t}_{0}}){C}'(\rho ,u(\rho ))C(\rho ,u(\rho ))\Phi (\rho ,{{t}_{0}}) \\
\times \left( \int_{{{t}_{0}}}^{\rho }{\Phi }({{t}_{0}},s)f\left( s,y(s),d(s),u(s) \right)\text{d}s \right)\text{d}\rho, t\in \left[ {{t}_{0}},\tau  \right]
\end{multline} 

We are in a position to provide our main assumptions together with the statement and proof of our main result.\\
\textbf{A1.} For system (1.1) we assume that there exists a nonempty subset $I$ of $\mathbb{R}_{\ge0}$ in such a way that for all ${{t}_{0}}\in I$, $\tau >{{t}_{0}}$ close to ${{t}_{0}}$ and for each $u\in {{U}_{[{{t}_{0}},\tau ]}}$ and $y\in {{Y}_{[{{t}_{0}},\tau ],u}}$ it holds: 
\begin{IEEEeqnarray}{rCl}
\Psi \left( t;{{t}_{0}},y(t),u(t)\right) >0,\forall t\in ({{t}_{0}},\tau ]\IEEEeqnarraynumspace \end{IEEEeqnarray} 
where the map $\Psi$ is given by (2.9); \\
\textbf{A2.} In addition, we assume that for every ${{t}_{0}}\in I$, $T>{{t}_{0}}$ close to ${{t}_{0}}$, $u\in {{U}_{[{{t}_{0}},T]}}$, $y\in {{Y}_{[{{t}_{0}},T],u}}$, $\ell \in (0,1)$ and constants $R, \theta >0$ there exists a constant $\tau \in \left( {{t}_{0}},\min \left\{ {{t}_{0}}+\theta ,T \right\}\right)$  such that 
\begin{multline}
\|{{\Psi }^{-1}}(\cdot ;{{t}_{0}},y,u)\left[ {{\Xi }}(\cdot ;{{t}_{0}},y,{{d}_{1}},u)-{{\Xi }}(\cdot ;{{t}_{0}},y,{{d}_{2}},u) \right]{{\|}_{({{t}_{0}},\tau ]}} \\ \le \ell \|{{d}_{1}}-{{d}_{2}}{{\|}_{[{{t}_{0}},\tau ]}}, \forall {{d}_{1}},{{d}_{2}}\in {{C}^{0}}\left( [{{t}_{0}},\tau ];{{\mathbb{R}}^{n}} \right),\\ \text{ with } \|{{d}_{i}}{{\|}_{[{{t}_{0}},\tau ]}}\le R,i=1,2
\end{multline} 

Assumption A1 is a type of persistence of excitation and A2 is a type of contraction condition. Assumptions A1 and A2 are in general difficult to be checked, however, they are both fulfilled for a class of nonlinear triangular systems, under weak assumptions that are exclusively expressed in terms of  system's dynamics (see (3.1) in the next section). We are in a position to state and prove our main result. Our approach leads to an explicit algorithm for the state estimation.

\begin{prop} Assume that A1 and A2 are fulfilled. Then the approximate SEDP is solvable for (1.1) over the set $I$; consequently (1.1) is observable over $I$.
\end{prop}
\begin{IEEEproof}[Proof of Proposition 2.1]
Let ${{t}_{0}}\in I$, $u\in {{U}_{[{{t}_{0}},\tau ]}}$, $y(\cdot )\in {{Y}_{[{{t}_{0}},\tau ],u}}$, with $\tau$ as given in A1 and A2, and let $x(\cdot ) \in C^0\left([t_0,\tau];\mathbb{R}^n\right) $ be a solution of (1.1) corresponding to $u(\cdot )$ and $y(\cdot )$ satisfying (2.3). Consider the trajectory ${{z}}(t)\text{:=}z(t;{{t}_{0}},{{z}}({{t}_{0}}),u)$ of the auxiliary system:
\begin{subequations}
\begin{IEEEeqnarray}{rCl}
\begin{aligned}
{{\dot{z}}}(t)=A\left( t,y(t),u(t) \right){{|}_{y(t)=C\left( t,u(t) \right)x(t)}}{{z}}(t) \IEEEeqnarraynumspace\IEEEeqnarraynumspace\\+f\left( t,y(t), {{z}}(t),u(t) \right){{|}_{y(t)=C\left( t,u(t) \right)x(t)}}\IEEEeqnarraynumspace\IEEEeqnarraynumspace
\end{aligned}\\
\text{with output }
{{Y}}(t):= C\left( t,u(t) \right){{z}}(t)\IEEEeqnarraynumspace\IEEEeqnarraynumspace\IEEEeqnarraynumspace
\end{IEEEeqnarray}
\end{subequations}
for certain initial $z(t_0) \in \mathbb{R}^n$. The map ${{Y}}(t)=C\left( t,u(t) \right){{z}}(t)$ is written: \begin{multline}
Y(t)= C\left( t,u(t) \right)\Phi (t,t_0)z(t_0)+C\left( t,u(t) \right)\\\times\int_{t_0}^{t}\Phi(t,s)f(s,y(s),z(s),u(s))\text{d}s \nonumber
\end{multline}
 By multiplying by ${\Phi }'(t,{{t}_{0}}){C}'\left( t,u(t) \right)$ and integrating we find:
\begin{multline*}
\int_{{{t}_{0}}}^{t}{{{\Phi }'}}(\rho ,{{t}_{0}}){C}'(\rho ,u(\rho )){{Y}}(\rho )\text{d}\rho \\= (\int_{{{t}_{0}}}^{t}{{{\Phi }'}}(\rho ,{{t}_{0}}){C}'(\rho ,u(\rho ))C(\rho ,u(\rho ))\Phi (\rho ,{{t}_{0}})\text{d}\rho \text{)}{{z}}({{t}_{0}})\\+\int_{{{t}_{0}}}^{t}({{{\Phi }'}}(\rho ,{{t}_{0}}){C}'(\rho ,u(\rho ))C(\rho ,u(\rho )) \Phi(\rho,t_0) \\ \times \int_{{{t}_{0}}}^{\rho} \Phi(t_0,s)f(s,y(s),z(s), u(s))\text{d}s)\text{d} \rho
\end{multline*}
The latter in conjunction with (2.9) - (2.11) yields:
\begin{multline}
{{z}}({{t}_{0}})={{\Psi }^{-1}}(t;{{t}_{0}},y,u)\int_{{{t}_{0}}}^{t}{{{\Phi }'}}(\rho ,{{t}_{0}}){C}'(\rho ,u(\rho )){{Y}}(\rho )\text{d}\rho \\ -{{\Psi }^{-1}}(t;{{t}_{0}},y,u){{\Xi }}(t;{{t}_{0}},y,z,u),\forall t\in ({{t}_{0}},\tau ] 
\end{multline}
By considering the solution $x(\cdot)$ of (1.1) corresponding to same $u(\cdot)$ with $y(t)=C(t,u(t))x(t), t \in [t_0,\tau]$ and with same initial $x(t_0)$, it follows from (2.13) that the mappings $Y(\cdot)$ and $y(\cdot)$ coincide, therefore, from (2.14) we get: 
\begin{multline}
{{x}}({{t}_{0}})={{\Psi }^{-1}}(t;{{t}_{0}},y,u)\int_{{{t}_{0}}}^{t}{{{\Phi }'}}(\rho ,{{t}_{0}}){C}'(\rho ,u(\rho )){{y}}(\rho )\text{d}\rho  \\ 
-{{\Psi }^{-1}}(t;{{t}_{0}},y,u){{\Xi }}(t;{{t}_{0}},y,x_{\varepsilon },u),\forall t\in ({{t}_{0}},\tau ]
\end{multline}
Let $T \in (t_0,\tau]$ and define:
\begin{multline}
\mathcal{F}_{T}(t;{{t}_{0}},y,z,u):=\Phi (t,{{t}_{0}}){{\Psi }^{-1}}(T;{{t}_{0}},y,u)\\ \times \left( \int_{{{t}_{0}}}^{T}{{{\Phi }'}}(\rho ,{{t}_{0}}){C}'(\rho ,u(\rho ))y(\rho )\text{d}\rho-\Xi(T;{{t}_{0}},y,z,u)\right) \\+\int_{{{t}_{0}}}^{t}{\Phi }(t,\rho )f(\rho ,y(\rho ), z(\rho ),u(\rho ))\text{d}\rho, \\
t\in [{{t}_{0}},T], z(\cdot )\in {{C}^{0}}\left( [{{t}_{0}},T];{{\mathbb{R}}^{n}} \right)
\end{multline}
Then, by (2.15) and (2.16) we have: 
\begin{equation}
{{\mathcal{F}}_{T}}(t;{{t}_{0}},y,{{x}},u)={{x}}(t),\forall t\in [{{t}_{0}},T]
\end{equation}
Next, consider a strictly increasing sequence $\left\lbrace R_{\nu} \in \mathbb{R}_{>0}, \nu=1,2,\ldots \right\rbrace$ defined as: 
\begin{equation}
R_{\nu+1}=2R_{\nu}, \nu=1,2,\ldots, \text{ wth } R_1=1
\end{equation}
Since, due to continuity of $x(\cdot)$, the set $\left\lbrace x(t), t \in [t_0,\tau] \right\rbrace$ is bounded, there exists an integer $k\ge 1$ such that
\begin{equation}
\Vert x(\cdot)\Vert_{[t_0,\tau]}<R_k
\end{equation}
Let 
\begin{equation}
\ell \in (0,1/2]
\end{equation}
By virtue of (2.1), (2.12) and (2.16) it follows that for the above $\ell$ there exists a decreasing continuous function $T=T(R):\mathbb{R}_{>0} \to (t_0,\tau]$ with $\lim_{R\to +\infty} T(R)=t_0$ and such that
\begin{multline}
\Vert\mathcal{F}_{T}(\cdot;{{t}_{0}},y,{{d}_{1}},u)-\mathcal{F}_{T}(\cdot ;{t_0},y,{{d}_{2}},u)\Vert_{(t_0,T ]} \le \ell \Vert d_{1}-d_2\Vert_{[{{t}_{0}},T ]}, \\T:=T(R), \forall d_1, d_2 \in C^{0}\left( [t_0, T(R)];\mathbb{R}^n\right), \text{ for which }\\ \max\left\lbrace \Vert d_i \Vert_{[t_0, T(R)]}, i=1,2 \right\rbrace \le R 
\end{multline}
Finally, define:
\begin{subequations}
\begin{IEEEeqnarray}{rCl}
t_{\nu}:=T(R_{\nu}); \IEEEeqnarraynumspace\IEEEeqnarraynumspace\\ 
\begin{aligned}
z_{\nu+1}(t):=\mathcal{F}_{t_{\nu}}(t;{{t}_{0}},y,{{z}_{\nu}},u);\IEEEeqnarraynumspace\\ t\in [{{t}_{0}},{{t}_{\nu}}], \nu=k,k+1,k+2,\ldots\IEEEeqnarraynumspace
\end{aligned}
\end{IEEEeqnarray}
\end{subequations}
with arbitrary
\begin{equation}
z_k \in C^0\left( [t_0,t_k];\mathbb{R}^n\right):\Vert z_k\Vert_{[t_0,t_k]}<R_k 
\end{equation}
Then by (2.19) and (2.21)-(2.23) we get
\begin{equation}
\Vert \mathcal{F}_{t_k}(\cdot;{{t}_{0}},y,z_k,u)-\mathcal{F}_{t_k}(\cdot ;{t_0},y,x,u)\Vert_{[t_0,t_k ]} \le \ell \Vert z_k-x\Vert_{[{{t}_{0}},t_k ]}
\end{equation}
According to (2.22b) let ${{z}_{k+1}}(\cdot ):={{\mathcal{F}}_{t_k}}(\cdot ;{{t}_{0}},y,{{z}_{k}},u)$. Then, from (2.17), (2.21)-(2.24) and the fact that the sequence $\left\{ {{t}_{\nu }}\in ({{t}_{0}},\tau] \right\}$ is decreasing, we have: 
\begin{multline}
\|{{z}_{k+1}}-x{{\|}_{[{{t}_{0}},{{t}_{k+1}}]}}=\|{{\mathcal{F}}_{{{{t}_{k}}}}}(\cdot ;{{t}_{0}},y,{{z}_{k}},u)-x{{\|}_{[{{t}_{0}},{{t}_{k+1}}]}}\\
=\|{{\mathcal{F}}_{{{{t}_{k}}}}}(\cdot ;{{t}_{0}},y,{{z}_{k}},u)-{{\mathcal{F}}_{{{t}_{k}}}}(\cdot ;{{t}_{0}},y,x,u){{\|}_{[{{t}_{0}},{{t}_{k+1}}]}}\\
\le \|{{\mathcal{F}}_{{{{t}_{k}}}}}(\cdot ;{{t}_{0}},y,{{z}_{k}},u)-{{\mathcal{F}}_{{{t}_{k}}}}(\cdot ;{{t}_{0}},y,x,u){{\|}_{[{{t}_{0}},{{t}_{k}}]}}\le \ell \|{{z}_{k}}-x{{\|}_{[{{t}_{0}},{{t}_{k}}]}}
\end{multline}
therefore, by invoking (2.18), (2.19), (2.20), (2.23) and (2.25) it follows:                                              
\begin{multline}
\|{{z}_{k+1}}{{\|}_{[{{t}_{0}},{{t}_{k+1}}]}}\le \|x{{\|}_{[{{t}_{0}},{{t}_{k+1}}]}}+\ell \|{{z}_{k}}-x{{\|}_{[{{t}_{0}},{{t}_{k}}]}}<{{R}_{k}}+2\ell {{R}_{k}}\\\le {{R}_{k+1}}\nonumber
\end{multline}
Quite similarly, by induction we get:
\begin{multline}
\|{{z}_{\nu +1}}-x{{\|}_{[{{t}_{0}},{{t}_{\nu +1}}]}}\le {{\ell }^{\nu +1-k}}\|{{z}_{k}}-x{{\|}_{[{{t}_{0}},{{t}_{1}}]}}; \\ 
\quad \|{{z}_{\nu }}{{\|}_{[{{t}_{0}},{{t}_{\nu +1}}]}}\le{{R}_{\nu }},\forall \nu =k,k+1,k+2,\ldots
\end{multline}
which implies
\begin{equation}
{{\xi }_{\nu }}:={{z}_{\nu}}({{t}_{\nu}})\to x({{t}_{0}});{{X}_{\nu }}({{t}_{0}},y,u):={{\xi }_{\nu }}
\end{equation}
where the values ${{\xi }_{\nu}}$ above are exclusively dependent on the values of the input $u$ and the output $y$ and the dynamics of system and they are independent of any time-derivatives of  $u$ and  $y$, thus both (2.6a) and (2.6b) are fulfilled. We conclude that the approximate SEDP is solvable for (1.1) over $I$, therefore, system (1.1) is observable over $I$.                                  
\end{IEEEproof}
The existence result of Proposition 2.1 does not in general determine explicitly the desired  sequence of mappings ${{X}_{\nu }}$ exhibiting (2.27). The reason is that, although existence of the constant $k$ satisfying (2.19) is guaranteed from boundedness of $\left\{ x(t),t{{\in }{[{{t}_{0}},\tau ]}} \right\}$, its precise determination requires knowledge of a bound of the previous set, which, in general, is not available. The rest part of this section is devoted for the establishment of a constructive algorithm, exhibiting the state determination. The corresponding procedure is based on the approach given in proof of Proposition 2.1 plus some appropriate modifications. 
\begin{center}
\textit{Algorithm}
\end{center}

To simplify the procedure, we distinguish two cases: \\
\textbf{Case I:} First, we assume that a bounded region of the state space is a priori known, where the unknown initial state of (1.1) belongs. Particularly, assume that for all ${{t}_{0}}\in I$, almost all $\tau >{{t}_{0}}$ near ${{t}_{0}}$ and input $u\in {{U}_{[{{t}_{0}},\tau ]}}$, an open ball ${{B}_{R}}$ of radius $R>0$ centered at zero is known,  such that the corresponding set of outputs of (1.1) is modified as follows:               
\begin{multline}
{{Y}_{[{{t}_{0}},\tau ],u}}:=\{y\in {{C}^{0}}\left( [{{t}_{0}},\tau ];{{\mathbb{R}}^{k}} \right):y(t)=C\left( t,u(t) \right) \\ 
\times x\left( t;{{t}_{0}},{{x}_{0}},u \right),\text{ a.e.}\ t\in [{{t}_{0}},\tau ],\text{ for certain }{{x}_{0}}\in {{B}_{R}}\}
\end{multline}
For the case above we adopt a slight modification of the approach used for the  proof of Proposition 2.1. Our proposed algorithm contains two steps:\\
\textbf{Step 1:} Define 
\begin{equation}
 R_1:=R
\end{equation}
where the latter is involved in (2.28). Notice that, due to the additional assumption (2.28), it follows that (2.19) holds with $k=1$ and for $\tau$
close to $t_0$. Next, consider a strictly increasing sequence $\left\{ {{R}_{\nu }}\in {{\mathbb{R}}_{>0}},\nu =1,2,\ldots \right\}$ satisfying the first equality of (2.18), namely, 
\begin{equation}
R_{\nu+1}=2R_\nu ,\nu =1,2,\ldots
\end{equation}
and with ${{R}_{1}}$ as above. We set $\ell =1/2$ and find a decreasing sequence ${{t}_{\nu }}\in {{\mathbb{R}}_{>0}},\nu =1,2,\ldots$ with ${{t}_{\nu }}\to {{t}_{0}}$ and  in such a way that
\begin{multline}
 \|{{\mathcal{F}}_{{{t}_{\nu }}}}(\cdot ;{{t}_{0}},y,{{d}_{1}},u)-{{\mathcal{F}}_{{{t}_{\nu }}}}(\cdot ;{{t}_{0}},y,{{d}_{2}},u){{\|}_{[{{t}_{0}},{{t}_{\nu }}]}}\le \ell \|{{d}_{1}}-{{d}_{2}}{{\|}_{[{{t}_{0}},{{t}_{\nu }}]}}, \\ 
 \quad \forall {{d}_{1}},{{d}_{2}}\in {{C}^{0}}\left( [{{t}_{0}},{{t}_{\nu }}];{{\mathbb{R}}^{n}} \right):\\\max \left\{ \|{{d}_{i}}{{\|}_{[{{t}_{0}},{{t}_{\nu }}]}},i=1,2 \right\}\le {{R}_{_{\nu }}},\nu =1,2,\ldots
\end{multline}\\
\textbf{Step 2:} Consider the sequence ${{z}_{\nu +1}}(t):={{\mathcal{F}}_{{{{t}_{\nu }}}}}(t;{{t}_{0}},y,{{z}_{\nu }},u),t\in [{{t}_{0}},t_{\nu}]$ with arbitrary initial ${{z}_{1}}(\cdot )\in {{C}^{0}}\left( [{{t}_{0}},{{t}_{k}}];{{\mathbb{R}}^{n}} \right)$ satisfying (2.23) with $k=1$ and set 
\begin{equation}
{{X}_{\nu }}({{t}_{0}},y,u):={{z}_{\nu }}({{t}_{\nu }}),\nu =1,2,\ldots
\end{equation} 
It then follows that (2.27) holds with  unique $x({{t}_{0}})\in {{B}_{R}}$ satisfying (2.3). Particularly, we  have:
\begin{multline}\|{{z}_{\nu +1}}-x{{\|}_{[{{t}_{0}},{{t}_{\nu +1}}]}}\le {{\ell }^{\nu }}\|{{z}_{1}}-x{{\|}_{[{{t}_{0}},{{t}_{1}}]}},\nu =1,2,\ldots \nonumber
\end{multline}
therefore, the sequence of mappings ${{X}_{\nu }}$, as defined by (2.32),  exhibits the state determination.  above  satisfies the desired (2.5) and (2.6).    \\
\textbf{Case II (General Case):}  We now provide an algorithm, which exhibits the state determination for the general case, without any additional assumption. The algorithm contains  two steps:\\
\textbf{Step 1}: Repeat the same procedure followed in Case I, with $R=1,2,3,\ldots$ and  construct a sequences of mappings 
\begin{subequations}
\begin{IEEEeqnarray}{rCl}
\begin{aligned}
 z_{\nu +1}^{i}(t):=\mathcal{F}_{t_{\nu}^i }(t;{{t}_{0}},y,{{z}_{\nu }},u),t\in [{{t}_{0}},t_{\nu}^i ],\IEEEeqnarraynumspace \\ 
z_{1}^{i}(\cdot ):=0,\nu =1,2,3,\ldots; i=1,2,3,\ldots\IEEEeqnarraynumspace
\end{aligned}
\end{IEEEeqnarray}
associated with appropriate decreasing sequences $\left\lbrace  t_{\nu}^i\in (t_0,\tau]\right\rbrace, i=1,2,3,\ldots$, with $\lim_{\nu\to\infty} t_{\nu}^i=t_0$, by pretending that $\|x(\cdot ){{\|}_{[{{t}_{0}},\tau ]}}<i,i=1,2,\ldots$ and in such a way that, if we define $\xi _{\nu }^{i}:=z_{\nu }^{i}(t_{\nu }^{i})$,  we have:
\begin{IEEEeqnarray}{rCl}
\begin{aligned}
 \left| \xi _{\nu +1}^{i}-x({{t}_{0}}) \right|\le {{\ell }^{\nu }}\|z_1^i-x{{\|}_{[{{t}_{0}},t_{1 }^{i} ]}}; \Vert z_1^i\Vert_{[t_0,t_{1}^{i}]}<i,\\ \forall \nu ,i=1,2,\ldots, 
\text{ provided that } \|x(\cdot ){{\|}_{[{{t}_{0}},\tau ]}}<i 
\end{aligned}
\end{IEEEeqnarray}
\end{subequations}
\textbf{Step 2}: Define 
\begin{equation}
X_{\nu}({{t}_{0}},y,u):=\xi _{\nu }^{\nu },\nu =1,2,\ldots
\end{equation}
Notice that, since the set $\left\{ x(t),t{{\in }{[{{t}_{0}},\tau ]}} \right\}$ is bounded, there exists an integer $k$ such that $\|x(\cdot ){{\|}_{[{{t}_{0}},\tau ]}}<k<k+1<k+2<\ldots$. The latter in conjunction with (2.33) yields: 
\begin{equation}\left| \xi _{\nu }^{\nu }-x({{t}_{0}}) \right|\le {{\ell }^{\nu-1 }}\|z_1^{\nu}-x{{\|}_{[{{t}_{0}},t^{\nu}_{1} ]}},\le 
\ell^{\nu-1}(\nu+k), \nu =k+1,k+2,\ldots\end{equation}
which, due to selection $\ell=1/2$, implies ${{X}_{\nu }}({{t}_{0}},y,u):=\xi _{\nu }^{\nu }\to x({{t}_{0}})$. We conclude that for the general case the sequence of mappings ${{X}_{\nu }}$, as defined by (2.34),  exhibits the state determination. Finally, it should be noted that,  according to the methodology above, contrary to the approach adopted in the proof of Proposition 2.1, the specific knowledge of $k$  satisfying (2.35) is not required.

\section{Application}
In this section we apply the results of Section II to triangular systems (1.1) of the form:
\begin{subequations}
\begin{IEEEeqnarray}{rCl}
\begin{aligned}
\dot{x}_i=a_{i+1}(t,x_1,u)x_{i+1}+f_i(t,x_1,u),i=1,\ldots,n-1,\\
\dot{x}_n=f_n(t,x_1,\ldots,x_n,u)\IEEEeqnarraynumspace\IEEEeqnarraynumspace\IEEEeqnarraynumspace\IEEEeqnarraynumspace
\end{aligned}
\\
(x_1,x_2,\ldots,x_n) \in \mathbb{R}^{n}, u \in \mathbb{R}^{m}, \IEEEeqnarraynumspace\IEEEeqnarraynumspace\IEEEeqnarraynumspace\nonumber
\end{IEEEeqnarray}
with output 
\begin{IEEEeqnarray}{rCl}
y=x_1
\end{IEEEeqnarray}
\end{subequations}
where we make the following assumptions:\\
\textbf{H1 (Regularity Assumptions).} It is assumed that for each $(x,y,u)\in {{\mathbb{R}}^{n}}\times \mathbb{R}\times {{\mathbb{R}}^{m}}$ the mappings ${{a}_{i}}(\cdot ,y,u),\,\,i=2,\ldots,n$, ${{f}_{i}}(\cdot ,{{x}_{1}},u)$, $i=1,\ldots ,n-1$ and ${{f}_{n}}(\cdot, {{x}_{1}},\ldots, {{x}_{n}},u)$ are measurable and locally essentially bounded and for each fixed $t \ge 0$ and $u\in {{\mathbb{R}}^{m}}$ the mappings ${{a}_{i}}(t,\cdot,u),i=2,\ldots ,n,{{f}_{i}}(t,\cdot ,u),i=1,\ldots ,n-1$ and ${{f}_{n}}(t,\ldots,u)$ are (locally) Lipschitz. 

Obviously, (3.1) has the form of (1.1) with 
\begin{subequations}
\begin{IEEEeqnarray}{rCl}
A(t,y,u):=\IEEEeqnarraynumspace\IEEEeqnarraynumspace\IEEEeqnarraynumspace\IEEEeqnarraynumspace\IEEEeqnarraynumspace\IEEEeqnarraynumspace\IEEEeqnarraynumspace\IEEEeqnarraynumspace\nonumber\\
{
\begin{bmatrix}
0&a_2(t,y,u)&0&\cdots&0\\
0&0&a_3(t,y,u)&\cdots&0\\
\vdots&\vdots&\ddots&\\
0&0&\cdots&&a_n(t,y,u)\\
0&0&0&\cdots&0
\end{bmatrix}
}\IEEEeqnarraynumspace\IEEEeqnarraynumspace\\
f:=[\begin{matrix}{{f}_{1}}&{{f}_{2}}&\cdots&{{f}_{n}}\end{matrix}]'\IEEEeqnarraynumspace\IEEEeqnarraynumspace\IEEEeqnarraynumspace\IEEEeqnarraynumspace \\
C\left( t,u \right):=C=\left[ \begin{matrix}
   1 & 0 & \cdots  & 0 
\end{matrix} \right]\IEEEeqnarraynumspace\IEEEeqnarraynumspace \IEEEeqnarraynumspace
\end{IEEEeqnarray}
\end{subequations}
We also make the following observability assumption:\\ 
\textbf{H2.} There exists a measurable set $I\subset {{\mathbb{R}}_{\ge 0}}$ with nonempty interior such that for all ${{t}_{0}}\in I$, $\tau >{{t}_{0}}$ close to ${{t}_{0}}$ and for each $u\in {{U}_{[{{t}_{0}},\tau ]}}:={{L}_{\infty }}([{{t}_{0}},\tau ];{{\mathbb{R}}^{m}})$ and $y \in Y_{[t_0,\tau],u}$ it holds: 
\begin{IEEEeqnarray}{rCl}
\prod\limits_{i=2}^{n}{{{a}_{i}}}({{t}_{0}})\ne 0; {{a}_{i}}({{t}_{0}})\text{ :=}{{a}_{i}}({{t}_{0}},y({{t}_{0}}),u({{t}_{0}}))\end{IEEEeqnarray}

\begin{prop} For the system (3.1) assume that H1 and H2 hold with ${{U}_{[{{t}_{0}},\tau ]}}:={{L}_{\infty }}([{{t}_{0}},\tau ];{{\mathbb{R}}^{m}})$ for certain $\tau >{{t}_{0}}$ close to ${{t}_{0}}\in I$. Then there exists a set $\hat{I}\subset I$ with $\text{cl}\hat{I}=I$ such that the approximate SEDP  is solvable over $\hat{I}$ for the system (3.1) by employing the methodology of Proposition 2.1; consequently, (3.1) is observable over $\hat{I}$.
\end{prop}
\begin{rem} A stronger version of assumption (3.3), is required in [8], [9], [25], [27], [28], for the construction of Luenberger type observers for a more general  class of triangular systems. Particularly, in all  previously mentioned works it is further imposed that the mappings $a_i(\cdot,\cdot,\cdot)$ are $C^1$. We note that the second conclusion of Proposition 3.1 concerning observability can alternatively be obtained under H1 and H2 as follows: By exploiting (3.3) and applying successive differentiation with respect to time, we can determine a map $X$ satisfying (2.4) with the information of the time-derivatives of $u$ and output $y$ (details are left to the reader). But this map is not acceptable for the solvability of SEDP for (3.1), due to  the additional requirements of Definition 2.2 that the candidate $X$ should be independent of the time-derivatives of $u$ and $y$. 
\end{rem}
\begin{IEEEproof}[Proof of Proposition 3.1]
We establish that the assumptions H1 and H2 guarantee that conditions A1 and A2 of previous section are fulfilled for (3.1), therefore by invoking Proposition 2.1 we get the desired statement. We first evaluate the fundamental solution $\Phi (t,{{t}_{0}})$ of (2.8) with $A(t,y,u)$ as given by (3.2a) for certain $t_0 \in I$. We find:
\begin{multline}
\Phi(t,t_0)=\\
{
\begin{bmatrix}
\varepsilon _{11}(=1)& \varepsilon_{12}(t)(t-t_0) & \cdots & \varepsilon_{1n}(t)(t-t_0)^{n-1}\\ 0 & \varepsilon_{21}(=1)& \cdots & \varepsilon_{2,n-1}(t)(t-t_0)^{n-2}\\ \vdots & \vdots & \ddots & \vdots \\ 0&0&\cdots& \varepsilon_{n1} (=1)
\end{bmatrix}
}
\end{multline}
where each function ${{\varepsilon }_{ij}}:[{{t}_{0}},\tau ]\to \mathbb{R}, i,j=1,2,\ldots,n$ satisfies 
\begin{subequations}
\begin{equation}
\varepsilon_{ij}(t)=E_{ij} \left( 1 +Z_{ij}(t) \right)
\end{equation}
for certain constants ${{E}_{ij}}\in \mathbb{R}$ and functions ${{Z}_{ij}}\in {{L}_{\infty }}([{{t}_{0}},\tau ];\mathbb{R})$ with 
\begin{equation}
\underset{t\to {{t}_{0}}}{\mathop{\lim }}\,{{Z}_{ij}}(t)=0
\end{equation}
\end{subequations}
Particularly, due to (3.3), we have:                            
\begin{IEEEeqnarray}{rCl}
\begin{aligned}
{{\varepsilon }_{1,1}}({{t}_{0}})=1, {{\varepsilon }_{1,i}}({{t}_{0}})={{E}_{1,i}}={{c}_{i}}{{a}_{2}}({{t}_{0}})\cdots {{a}_{i}}({{t}_{0}})\ne 0,\IEEEeqnarraynumspace\\ i=2,\ldots,n \IEEEeqnarraynumspace\IEEEeqnarraynumspace\IEEEeqnarraynumspace\IEEEeqnarraynumspace\IEEEeqnarraynumspace\IEEEeqnarraynumspace\IEEEeqnarraynumspace\IEEEeqnarraynumspace\IEEEeqnarraynumspace
\end{aligned}
\end{IEEEeqnarray}
for certain nonzero constants ${{c}_{i}}$. Notice, that the above representation is feasible for almost all $t_0 \in I$ due to our regularity assumptions concerning ${{a}_{i}}$. For simplicity, we may assume next that (3.4) - (3.6) hold for every $t_0 \in I$. We now calculate by taking into account (3.2c) and (3.4):
\begin{IEEEeqnarray}{rCl}
C\Phi \left( t,{{t}_{0}} \right)=[{{\varepsilon }_{11}}(t),{{\varepsilon }_{12}}(t)(t-{{t}_{0}}),\cdots,{{\varepsilon }_{1n}}(t){{(t-{{t}_{0}})}^{n-1}}]\IEEEeqnarraynumspace
\end{IEEEeqnarray}
Notice that $\Psi$, as defined by (2.9), satisfies (2.11) since otherwise, there would exist sequences $\upsilon^i=(\upsilon_1^i,\upsilon_2^i,\ldots,\upsilon_n^i) \in \mathbb{R}^n\setminus\lbrace0\rbrace, { {{t}_{i}}\in ({{t}_{0}},\tau ]}$ and a nonzero vector $\upsilon =[{{\upsilon }_{1}},{{\upsilon }_{2}},\cdots,{{\upsilon }_{n}}]\in {{\mathbb{R}}^{n}}$ with $\lim_{i\to \infty} \upsilon^i=\upsilon, \lim_{i\to \infty} {{t}_{i}}={{t}_{0}}$ and in such a way that $C\Phi ({{t}_{i}},{{t}_{0}})\upsilon^i =0,i=1,2,\ldots$. Then by using (3.7) we get ${{\upsilon }_{1}^i}{{\varepsilon }_{11}}({{t}_{i}})+{{\upsilon }_{2}^i}{{\varepsilon }_{12}}({{t}_{i}})({{t}_{i}}-{{t}_{0}})+\ldots+{{\upsilon }_{n}^i}{{\varepsilon }_{1n}}({{t}_{i}}){{({{t}_{i}}-{{t}_{0}})}^{n-1}}=0$, which by virtue of (3.5) and (3.6) implies that $\upsilon =0$, a contradiction. We conclude that relation (2.11) of A1 holds with ${{U}_{[{{t}_{0}},\tau ]}}={{L}_{\infty }}([{{t}_{0}},\tau ];{{\mathbb{R}}^{m}})$.  In order to establish A2, we calculate, according to definition (2.9) and by using (3.7):
\begin{equation}
\Psi=
{
\begin{bmatrix}
\varepsilon_{11}(t)\Delta t&\varepsilon_{12}(t)(\Delta t)^2&\cdots & \varepsilon_{1n}(t)(\Delta t)^n\\
\varepsilon_{21}(t)(\Delta t)^2&\varepsilon_{22}(t)(\Delta t)^3&\cdots & \varepsilon_{2n}(t)(\Delta t)^{n+1}\\
\vdots&\vdots&&\vdots\\
\varepsilon_{n1}(t)(\Delta t)^{n}&\varepsilon_{n2}(t)(\Delta t)^{n+1}&\cdots & \varepsilon_{nn}(t)(\Delta t)^{2n-1}
\end{bmatrix}
}
\end{equation}
where $\Delta t:=t-t_0$ and the functions ${{\varepsilon }_{ij}}:[{{t}_{0}},\tau ]\to \mathbb{R}$ above satisfy (3.5). Define ${{d}_{1}}=({{x}_{2}},\ldots ,{{x}_{n}}{)}'$, ${{d}_{2}}=({{\bar{x}}_{2}},\ldots ,{{\bar{x}}_{n}}{)}'$ and let $f:=[f_1,f_2,\cdots,f_{n-1},f_n]'$ and 
\begin{IEEEeqnarray}{rCl}
\Delta f(\cdot ,{{d}_{1}},{{d}_{2}},u):=f(\cdot ,{{d}_{2}},u)-f(\cdot ,{{d}_{1}},u)
\end{IEEEeqnarray}
By (3.1a), (3.2c), (3.4) and (3.9) we find:
\begin{multline}
\int_{{{t}_{0}}}^{\rho }{\Phi }(t_0,s)\Delta f(s,{{d}_{1}}(s),{{d}_{2}}(s),u(s))\text{d}s =\\ \left[ \varepsilon_1(\rho)(\Delta \rho)^n,\varepsilon_2(\rho)(\Delta \rho)^{n-1},\cdots,\varepsilon_{n-1}(\rho)(\Delta \rho)^2, \varepsilon_n(\rho)\Delta \rho \right]'
\end{multline}
$\rho \ge {{t}_{0}}$ near ${{t}_{0}}$, $\Delta \rho :=t_0-\rho$, where the functions ${{\varepsilon }_{i}}(\cdot )={{\varepsilon }_{i}}(\cdot ;{{d}_{1}}(\cdot ),{{d}_{2}}(\cdot ))$ have  the form: 
\begin{IEEEeqnarray}{rCl}
{{\varepsilon }_{i}}(t)={{E}_{i}}(1+{{Z}_{i}}(t)),\ li{{m}_{t\to {{t}_{0}}}}{{Z}_{i}}(t)=0, i=1,\ldots,n \IEEEeqnarraynumspace
\end{IEEEeqnarray}
for certain ${{E}_{i}}\in \mathbb{R},$ ${{Z}_{i}}\in {{L}_{\infty }}([{{t}_{0}},\tau ];\mathbb{R})$, and in such a way that, due to (3.9) and Lipschitz continuity property of ${{f}_{n}}$, the following holds for every $R>0$:  
\begin{multline}
\vert\varepsilon_i(t)\vert \le C \Vert d_1-d_2 \Vert_{[t_0,\tau]},  \forall t \in [t_0,\tau], \tau \text{ near } t_0, \\ {{d}_{1}},{{d}_{2}}\in {{C}^{0}}\left( [{{t}_{0}},\tau ];{{\mathbb{R}}^{n}} \right) \text{ with } \|{{d}_{i}}{{\|}_{[{{t}_{0}},\tau ]}}\le R,i=1,2 
\end{multline}
for certain constant $C>0$. Also, we evaluate from (3.8): 
\begin{multline}
\Psi^{-1}=\\
\begin{bmatrix}
\varepsilon_{11}(t)(\Delta t)^{-1}&\varepsilon_{12}(t)(\Delta t)^{-2}&\cdots&\varepsilon_{1n}(t)(\Delta t)^{-n}\\
\varepsilon_{21}(t)(\Delta t)^{-2}&\varepsilon_{22}(t)(\Delta t)^{-3}&\cdots&\varepsilon_{2n}(t)(\Delta t)^{-n-1}\\
\vdots&\vdots&&\vdots\\
\varepsilon_{n1}(t)(\Delta t)^{-n}&\varepsilon_{n2}(t)(\Delta t)^{-n-1}&\cdots &\varepsilon_{nn}(t)(\Delta t)^{-2n+1}
 \end{bmatrix}
\end{multline}
where $\Delta t:=t-t_0$ and ${{\varepsilon }_{ij}}$ above satisfy (3.5a). From (2.10), (3.7) and (3.10) - (3.12) we also find: 
\begin{multline}
\Xi (t;{{t}_{0}},y,{{d}_{1}},u)-\Xi (t;{{t}_{0}},y,{{d}_{2}},u)= \\ 
{{\left[ {{\varepsilon }_{1}}(t){{(\Delta t)}^{n+1}},{{\varepsilon }_{2}}(t){{(\Delta t)}^{n+2}},\cdots,{{\varepsilon }_{n}}(t){{(\Delta t)}^{2n}} \right]}^{\prime }}  
\end{multline}
for $t$ near ${{t}_{0}}$, where $\Delta t:=t-t_0$ and each ${{\varepsilon }_{i}},i=1,\ldots ,n$ above satisfy again (3.11) and (3.12). The latter in conjunction with (3.9), (3.13) and (3.14) implies A2. To be precise, the following holds: For every ${{t}_{0}}\in I$, $T>{{t}_{0}}$ close to ${{t}_{0}}$, $u\in {{U}_{[{{t}_{0}},T]}}$, $y\in {{Y}_{[{{t}_{0}},T],u}}$ and constants $\ell \in (0,1)$ and $R,\theta >0$, a constant $\tau \in ({{t}_{0}},\min \left\{ {{t}_{0}}+\theta ,T \right\})$ can be found satisfying (2.12).  We conclude that both A1 and A2 are fulfilled for the case (3.1), hence, according to Proposition 2.1, the approximate SEDP is solvable for (3.1) over a set $\hat{I}\subset I$ with $\text{cl}\hat{I}=I$.
\end{IEEEproof}
\begin{exm} We illustrate the nature of our methodology by considering the elementary case of the planar single-input triangular system ${{{\dot{x}}}_{1}}={{x}_{2}}u,  
{{{\dot{x}}}_{2}}={{x}_{1}}-x_{2}^{3}$ with output $y={{x}_{1}}$ that has the form (3.1) with $$A:=\left[ \begin{matrix}
   0 & u  \\
   0 & 0  \\
\end{matrix} \right] \text{and } f:=\left( 0,{{x}_{1}}-x_{2}^{3} \right)'.$$ We may assume that each admissible input $u$ is any nonzero measurable and essentially locally bounded function and for simplicity, let $u(t)=1$ for $t$ near zero. Obviously, the system above satisfies H1 and H2. Let us choose $({{x}_{1}}(0),{{x}_{2}}(0))=(2,\,0)$ as initial condition and calculate the corresponding output trajectory $y={{x}_{1}}$ (see Figure 1 below).
\begin{figure}[h]
\caption{Output of the System}
\centering
\includegraphics[scale=0.6]{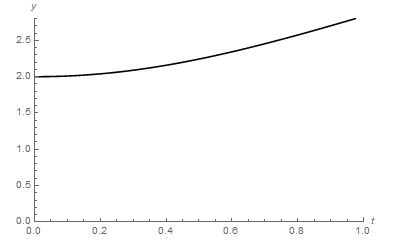}
\end{figure}
We next apply the methodology suggested in the previous section,  in order to confirm that our proposed algorithm converges to $x(0)=({{x}_{1}}(0),{{x}_{2}}(0))$ above. For simplicity, let us assume that is a priori known that the “unknown" initial state $x(0)$ is contained into the ball ${{B}_{R}}$ of radius $R=3$ centered at zero. We take ${{R}_{1}}=3,{{R}_{2}}=6,{{R}_{3}}=12,{{R}_{4}}=24,{{R}_{5}}=48,\ldots $ and $\ell=0.5$ as in the proposed algorithm (Case I). By taking into account the known values of $y(\cdot )$, we find a decreasing sequence $\left\{ {{t}_{\nu }} \right\}$ satisfying  (2.31)  converging to $t_0=0$; particularly, take ${{t}_{1}}=5\times {{10}^{-4}}$ ${{t}_{2}}=1.3\times {{10}^{-4}},\ {{t}_{3}}=3.2\times {{10}^{-5}},\ {{t}_{4}}=7.8\times {{10}^{-6}},\ {{t}_{5}}=1.9\times {{10}^{-6}},\ldots $. Then, choose an arbitrary constant initial map, say ${{z}_{1}}(\cdot )=(z_{1}^{1}(\cdot ),z_{1}^{2}(\cdot ))$; with $z_{1}^{1}\left( t \right)=0$ and $z_{1}^{2}\left( t \right)=1$, for $t\in [{{t}_{0}}=0,{{t}_{1}}=5\times {{10}^{-4}}]$  and successively apply (2.22b), in order to evaluate the desired sequence  $\left\{ {{z}_{\nu }}({{t}_{\nu }}):=(z_{\nu }^{1}({{t}_{\nu }}),z_{\nu }^{2}({{t}_{\nu }}) \right\},\nu =1,2,\ldots $.  Figures 2 and 3 present the corresponding values of errors $e_{\nu }^{1}({{t}_{\nu }}):=z_{\nu }^{1}({{t}_{\nu }})-{{x}_{1}}(0),e_{\nu }^{2}({{t}_{\nu }}):=z_{\nu }^{2}({{t}_{\nu }})-{{x}_{2}}(0),\nu =1,2,\ldots$ and confirms that the evaluated pair of terms ${{z}_{\nu }}\left( {{t}_{\nu }} \right)$ converges to the pair  $({{x}_{1}}(0),{{x}_{2}}(0))=(2,0)$. 
\begin{figure}[h]
\caption{Error $e_{\nu}^1$}
\centering
\includegraphics[scale=0.6]{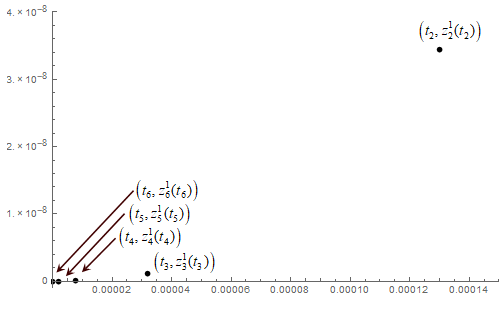}
\end{figure}
\begin{figure}[h]
\caption{Error $e_{\nu}^2$}
\centering
\includegraphics[scale=0.6]{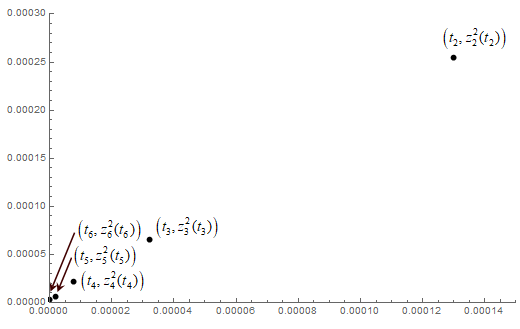}
\end{figure}

Finally, we remark that, since the system is forward complete, then for any $T>{{t}_{0}}$ the sequence of mappings ${{X}_{\nu }}(0,y,u):=x(t;0,{{\xi }_{\nu }},u),t\in [0,T],\nu =1,2,\ldots$ with ${{\xi }_{\nu }}:={{z}_{\nu }}\left( {{t}_{\nu }} \right)$ uniformly approximates the unknown solution $x(\cdot ;0,x(0),u)$ on  the interval $[{{t}_{0}},T]$.
\end{exm}
\section{Additional Hypotheses and Hybrid Observer}
In this section we briefly present a hybrid-observer technique for the state estimation for (1.1). The proof of the following proposition is based on a modification of the approach employed in Section II.
\begin{prop}
For the system (1.1) we make the same assumptions with  those imposed in statement of Proposition 2.1. Also, assume that for any ${{t}_{0}}\in I$ and  input  $u\in {{U}_{[{{t}_{0}},+\infty )}}$ there exists a constant $C>0$ such that
\begin{IEEEeqnarray}{rCl}
\left| f\left( t,y,{{z}_{1}},u(t) \right)-f\left( t,y,{{z}_{2}},u(t) \right) \right|\le C\left| {{z}_{1}}-{{z}_{2}} \right|,\nonumber \IEEEeqnarraynumspace\\ \forall ( t,y) \in [t_0,\infty) \times \mathbb{R}, {{z}_{1}},{{z}_{2}}\in \mathbb{R}^n \IEEEeqnarraynumspace
\end{IEEEeqnarray} 
Then, there exists a sequence $\xi_{\nu}=\xi_{\nu}(t_0,y,u) \in \mathbb{R}^n, \nu=0,1,2,\ldots$ such that, if for any arbitrary constant $h>t_0$ we define: 
\begin{subequations}
\begin{IEEEeqnarray}{rCl}
\begin{aligned}
{{\omega }_{0}}(\xi ):=\xi; \IEEEeqnarraynumspace\IEEEeqnarraynumspace\IEEEeqnarraynumspace\IEEEeqnarraynumspace\\ 
{{\omega }_{\nu+1}}(\xi )=x({{t}_{0}}+(\nu+1)\sigma;{{t}_{0}}+\nu\sigma ,{{\omega }_{\nu}}(\xi ),u),\nu=0,1,2,\ldots 
\end{aligned}\\
{{m}_{0}}:={{\omega }_{0}}({{\xi }_{0}})(={{\xi }_{0}});\ {{m}_{\nu}}:={{\omega }_{\nu}}({{\xi }_{\nu}}),\nu=1,2,\ldots\IEEEeqnarraynumspace
\end{IEEEeqnarray}
\end{subequations}
where $\sigma:=h-t_0$, then the system below exhibits the global state estimation of (1.1):
\begin{subequations}
\begin{equation}
\dot{\hat{x}}\left( t \right)=A(t,y,u)\hat{x}+f(t,y,\hat{x},u), \ t\in \left[ {{t}_{0}}+\nu\sigma ,{{t}_{0}}+(\nu+1)\sigma  \right)
\end{equation}
\begin{equation}
\hat{x}\left( {{t}_{0}}+\nu\sigma  \right)={{m}_{\nu}},\ \nu=0,1,2,\ldots
\end{equation}
\end{subequations}
particularly, it holds: 
\begin{IEEEeqnarray}{rCl}
\underset{t\to \infty }{\mathop{\lim }}\vert\hat{x}(t;t_0,m_0,u)-x(t;t_0,x_0,u)\vert=0\IEEEeqnarraynumspace
\end{IEEEeqnarray}  
\end{prop}
\begin{IEEEproof}[Outline of Proof]
Let $x(\cdot )=x(\cdot ;{{t}_{0}},{{x}_{0}},u),{{t}_{0}}\in I$ be a solution of (1.1) corresponding to $u(\cdot)$. Let  $h>{{t}_{0}}$ and consider the sequence ${{C}_{\nu }}:={{(\exp hC)}^{\nu +1}},\nu =1,2,\ldots$, with $C$ as defined in (4.1), and let $\left\{ {{\ell }_{\nu }} \right\},{{\ell }_{\nu }}\in (0,1/2],\nu =1,2,\ldots$ be a decreasing sequence with  
\begin{IEEEeqnarray}{rCl}
\lim_{\nu \to \infty}\ell_{\nu-1}C_{\nu}=0
\end{IEEEeqnarray} 
We next proceed by using a generalization of the procedure employed for the proof of Proposition 2.1. First, we find a decreasing sequence ${{t}_{\nu }}\in ({{t}_{0}},h]$, $,\nu =0,1,2,\ldots $ with ${{t}_{\nu }}\to {{t}_{0}}$ such that 
\begin{multline}
\|{{\mathcal{F}}_{{{t}_{\nu }}}}(\cdot ;{{t}_{0}},y,{{d}_{1}},u)-{{\mathcal{F}}_{{{t}_{\nu }}}}(\cdot ;{{t}_{0}},y,{{d}_{2}},u){{\|}_{[{{t}_{0}},{{t}_{\nu }}]}}\le {{\ell }_{\nu }}\|{{d}_{1}}-{{d}_{2}}{{\|}_{[{{t}_{0}},{{t}_{\nu }}]}}, \\ 
 \forall {{d}_{1}},{{d}_{2}}\in {{C}^{0}}\left( [{{t}_{0}},{{t}_{\nu }}];{{\mathbb{R}}^{n}} \right):\max \left\{ \|{{d}_{i}}{{\|}_{[{{t}_{0}},{{t}_{\nu }}]}},i=1,2 \right\}\le {{R}_{\nu }},\\\nu =k,k+1,k+2,\ldots  
\end{multline}
with $\left\{ {{R}_{\nu }} \right\}$ and $k$ as  defined by (2.18) and (2.19), respectively.  Then, consider the sequence of mappings ${{z}_{\nu }}(\cdot )\in {{\mathbb{R}}^{n}},\nu =k,k+1,k+2,\ldots $, as precisely defined in  (2.22) and again define: 
\begin{equation}
{{\xi }_{\nu }}:={{z}_{\nu }}({{t}_{\nu }}),\nu =k,k+1,\ldots
\end{equation}
Then, as in proof of Propostion 2.1 we can show, by exploiting (4.6), that 
\begin{equation}
\left| {{\xi }_{\nu }}-{{x}_{0}} \right|\le {{\ell }_{\nu }}{{\ell }_{\nu -1}}\cdots{{\ell }_{k}}\left| {{\xi }_{k}}-{{x}_{0}} \right|,\forall \nu =k+1,k+2,\ldots
\end{equation}
We are now in a position to show (4.4). We take into account (4.1)-(4.3), (4.5), (4.8), definition of ${{C}_{\nu }}$  and consider  the difference between the integral representation of the solutions of (1.1a) and (4.3a). Then, by successively applying  the  Gronwall - Bellman inequality, we can estimate:  
\begin{multline}
{{\left\| \hat{x}(\cdot ;{{t}_{0}}+\nu \sigma,{{m}_{\nu }},u)-x(\cdot ;{{t}_{0}},{{x}_{0}},u) \right\|}_{[{{t}_{0}}+\nu \sigma,{{t}_{0}}+(\nu +1)\sigma]}} \\ 
\le {{\ell }_{k}}{{\ell }_{k+1}}\cdots{{\ell }_{\nu -1}}{{C}_{\nu }}\left| {{\xi }_{0}}-{{x}_{0}} \right|\le {{\ell }_{\nu -1}}{{C}_{\nu }}\left| {{\xi }_{0}}-{{x}_{0}} \right|,\\\nu =k+1,k+2,\ldots 
\end{multline}
and the above, in conjunction with (4.5), asserts that                             \[\underset{\nu \to \infty }{\mathop{\lim }}\,{{\left\| \hat{x}(\cdot ;{{t}_{0}}+\nu \sigma,{{m}_{\nu }},u)-x(\cdot ;{{t}_{0}},{{x}_{0}},u) \right\|}_{[{{t}_{0}}+\nu \sigma,{{t}_{0}}+(\nu +1)\sigma]}}=0\] The latter implies the desired (4.4). Details are left to the reader. 
\end{IEEEproof}
\nocite{*}
\section{Conclusion}
Sufficient conditions for observability and solvability of the state estimation for a class of nonlinear control time - varying systems are derived. The state estimation is exhibited by means of a sequence of functionals approximating the unknown state of the system on a given bounded time interval. Each functional is exclusively dependent on the dynamics of system,  the input $u$ and the corresponding output $y$. The possibility of solvability of the state estimation problem by means of hybrid observers is briefly examined.

\bibliographystyle{IEEEtran}

\end{document}